\documentclass{article}
\usepackage[a4paper,left=2.5cm,right=2cm,top=2.5cm,bottom=2cm]{geometry}
\usepackage{graphicx,amsmath,latexsym,amssymb,amsthm,multirow}
\usepackage[colorlinks=black,urlcolor=black]{hyperref}
\usepackage{tfrupee}
\usepackage{enumitem}

\begin{document}

\begin{center}
\textbf{LEAKAGE INVENTORY MODEL WITHOUT SHORTAGES UNDER FUZZY PARAMETERS}\\
\vskip 0.5cm
HUIDROM MALEMNGANBI\\
Department of Mathematics, Manipur University, Canchipur-795003, India\\
email: malem\_feb82012@yahoo.in\\
\vskip 0.5cm
M.KUBER SINGH\\
Department of Mathematics, D M College of Science, Imphal-795001, India\\
email: moirang1@yahoo.com\\
\end{center}
\noindent {\bf Abstract:} In this paper, an attempt has been made to develop a simple leakage inventory model without shortages with instantaneous or finite production rate under fuzzy environment. In the present day scenario, it is difficult to decide the exact annual demand rate and hence a major difficulty is faced by a decision maker to forecast the demand. Also, in any inventory system, goods in stock may subject to deterioration or leakage. Deterioration of goods refers to decrease in quality whereas quantity remains more or less the same. On the other hand, leakage refers to loss in quantity whereas the quality remains unchanged for a certain period of time. Leakages in the inventory system may be considered to be very small, not detectable by the management immediately and hence it is difficult to decide the exact leakage rate. The objective of this paper is to consider these variable parameters and determine the optimal economic order quantity (EOQ) to maximise the annual total profit. So, fuzzy inventory models have been proposed considering fuzzy leakage rate and fuzzy annual demand rate to estimate the total profit per unit time. Signed distance method is used for defuzzification. Numerical examples are provided to support the results of the proposed models.
\vskip 0.2cm
\noindent {\bf Keywords:}  Leakage, Deterioration, EOQ, Fuzzy inventory model, Defuzzification\\

\noindent {\bf 2010 AMS Subject Classification:} 90B05

\section{Introduction}
In an inventory system dealing with liquid or gaseous stock, leakage is a realistic and common phenomenon. Leakage may be defined as the process in which material(liquid or gas) is lost through a leak in the storage facility or in-transit. This reduces the quantity stock and hence the loss due to this leakage can't be ignored while developing or discussing the economic order quantity(EOQ) or production order quantity(POQ) inventory models. In literature, many inventory models referring to deterioration under different conditions and circumstances have been developed. As an extension of inventory models of deteriorating items, numerous works on perishable items, evaporating items, imperfect items, etc. have been formulated. After all, deterioration, perishable, evaporation, leakage etc. all refer to the loss of the value of the stock either qualitatively or quantitatively. Taleizadeh\cite{talei} considered an EOQ model with partial backordering and advance payments for an evaporating item. So, considering leakage as analogous to deterioration, Tomba and Geeta \cite{tomba} developed a leakage inventory model having no shortages with uniform demand rate and instantaneous production rate.\\

In developing the models, there are many asumptions and parameters which are considered to be fixed and has exact value. But in reality, these assumptions seem to be unrealistic since they are generally vague and imprecise, sometimes even impossible to determine or guarantee the exact value. In this context, fuzzy inventory models are discussed as an extension of the traditional or crisp inventory models.  Yao and Lee\cite{yaolee} solved the inventory model with shortages by fuzzifying the order quantity using extension principle. Lee and Yao\cite{lee} also discussed production inventory problems by defuzzyfying demand quantity and production quantity considering triangular fuzzy numbers. Chang\cite{schang} developed a fuzzy production inventory for fuzzy product quantity using triangular fuzzy number. Yao and Su\cite{yaosu} studied the fuzzy inventory model with backorder for fuzzy total demand based on interval-valued fuzzy set. Chang\cite{chang} discussed the EOQ model with imperfect quality items by applying the set theory considering fuzzy defective rate and fuzzy annual demand. Mahata\cite{mahata} applied fuzzy set theory in an EOQ model for items with imperfect quality and shortage backordering. Jaggi et al.\cite{jaggi} studied a fuzzy inventory model for deteriorating items with initial inspection and allowable shortage under the condition of permissible delay in payments. Patro et al.\cite{patro} in their paper discussed an EOQ model for fuzzy defective rate with allowable proportionate discount. Chen and Ouyang\cite{chen}, in their paper extended an ordering policy for deteriorating items with allowable shortage and permissible delay in payment to fuzzy model by fuzzifying the carrying cost, interest paid rate and interest earned rate simultaneously based on the interval-valued fuzzy numbers and triangular fuzzy numbers. Chang et al.\cite{changetal} discussed a fuzzy mixed inventory model invoving variable lead-time with backorder and lost sales using probabilistic fuzzy set and triangular fuzzy number. They used two methods of defuzzification.\\

In this study, we investigate the leakage inventory model of Tomba and Geeta\cite{tomba} and proposed a fuzzy leakage model by considering demand rate and leakage rate as triangular fuzzy numbers. Yao and Wu's\cite{yaowu} ranking method for fuzzy numbers is used to find the optimal order quantity and the minimum total cost per unit time in fuzzy environment.  
\section{Brief review of Tomba and Geeta's model}
Tomba and Geeta \cite{tomba} developed a deterministic leakage inventory model without shortages having uniform demand rate with instantaneous production. They assumed the lead time of an order quantity q to be zero in developing their model.
\begin{center}
\includegraphics[width=0.4\textwidth]{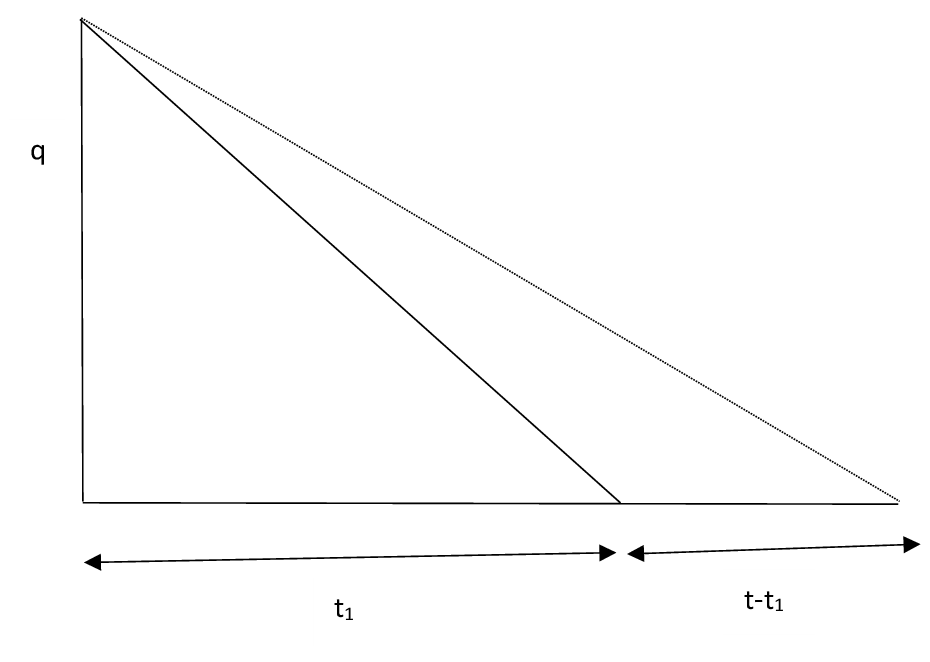}\\
Figure 1: Instantaneous production leakage inventory without shortage\\
\end{center}
According to their model, total holding cost = $\frac{1}{2}htq$, where h and q denote the holding cost per unit quantity per unit time t and order quantity per run respectively.\\
Average total cost per unit time = $\frac{1}{t}$(holding cost+set-up cost)=$\frac{1}{2}hq+\frac{s\phi}{q}$, where s is the set-up cost per order and $\phi$ is the demand rate per unit time respectively.\\
The model has a leakage with a leakage rate $\psi$ per unit time. So, the order quantity q reduces to zero in time $t_1$ and $q=(\phi+\psi)t_1$.\\
Then, $\phi t=q=(\phi+\psi)t_1 \implies dt-(\phi+\psi)t_1 = 0 \implies t-t_1=\frac{\psi q}{\phi(\phi+\psi)}$  \\
Quantity lost due to leakage=$ (t -t_1)\phi=\frac{\psi q}{(\phi+\psi)}$\\
So, average additional cost incurred due to leakage =  $\frac{(t-t_1)hq}{2t} =\frac{\psi hq}{2(\phi+\psi)}$\\
$\therefore$ Average total cost per unit time, 
\begin{align}
TC=(\frac{qh}{2})+(\frac{\phi s}{q})+\frac{\psi qh}{2(\phi+\psi)}
\end{align}
\section{Definitions and Preliminaries}
To extend the Tomba and Geeta's model to fuzzy environment, the following definitions and preliminaries are taken into account:
\subsection{Fuzzy point}
A fuzzy set $\tilde{b}$ defined on $\mathbb{R}$ is called a fuzzy point if its membership function is given by
\begin{center}
 $\mu_{\tilde{b}}(y) =\begin {cases} 1 & \mbox {if~~} y=b \\ 0 & \mbox {if~~} y \neq b \end{cases}$
 \end{center}
\subsection{$\alpha$-level fuzzy point}
Let $\tilde{b}_{\alpha}$ ($0\leqslant\alpha\leqslant1$) be a fuzzy set defined on $\mathbb {R}$. Then, it is called an $\alpha$-level fuzzy point if its membership function is defined as\\
\begin{center}
 $\mu_{{\tilde{b}}_{\alpha}}(y) =\begin {cases} \alpha & \mbox {if~~} y=b \\ 0 & \mbox {if~~} y \neq b \end{cases}$
 \end{center}
\subsection{Triangular Fuzzy Number}
A fuzzy set $\tilde{B} = (\beta_1,\beta_2,\beta_3)$ defined on $\mathbb{R}$, where $\beta_1<\beta_2<\beta_3$ is called a triangular fuzzy number if its membership function is given by
\begin{center}
$\mu_{\tilde{B}} (y)= \begin{cases}0 & \mbox {if  } y<\beta_1 \\ \frac {y-\beta_1}{\beta_2-\beta_1} & \mbox {if  } \beta_1 \leqslant y \leqslant \beta_2\\ \frac {\beta_3-y}{\beta_3-\beta_2} & \mbox {if  } \beta_2 \leqslant y \leqslant \beta_3\\ 0 & \mbox {if  }   y > \beta_3 \end{cases}$\\
\end{center}
\subsection{$\alpha$-level fuzzy interval}
For $0\leqslant\alpha\leqslant1$, the fuzzy set [$c_{\alpha},d_{\alpha}$] defined on $\mathbb{R}$ is called an $\alpha$-level fuzzy interval if the membership function of [$c_{\alpha},d_{\alpha}$] is given by \\
\begin{center}
$\mu_{[c_{\alpha},d_{\alpha}]}(y) =\begin {cases} \alpha & \mbox {if~~} c\leqslant y\leqslant d \\ 0 & \mbox {if~~}$ otherwise $\end{cases}$
\end{center}
\subsection{$\alpha$-cut of a fuzzy set}
Let $\tilde{B}$ be a fuzzy set on $\mathbb{R}$ and $0\leqslant\alpha\leqslant1$. Then, the $\alpha$-cut $B(\alpha)$ of $\tilde{B}$ consists of points y such that $\mu_{\tilde{B}}(y)\geqslant\alpha$, i.e., $B(\alpha)=\{y|\mu_{\tilde{B}}(y)\geqslant\alpha\}$.
\subsection{Decomposition Principle \cite{yaowu}} 
Let $\tilde{B}$ be a fuzzy set on $\mathbb{R}$ and $0\leqslant\alpha\leqslant1$. Let us suppose the $\alpha$-cut of $\tilde{B}$ as a closed interval $[B_L(\alpha),B_U(\alpha)]$, i.e., $B(\alpha)=[B_L(\alpha),B_U(\alpha)]$. Then,
\begin{align}\tilde{B}=\bigcup\limits_{0\leqslant\alpha\leqslant1}\alpha B(\alpha)=\bigcup\limits_{0\leqslant\alpha\leqslant1}[B_L(\alpha)_{\alpha},B_U(\alpha)_{\alpha}] 
\end{align}
or
\begin{align}
 \mu_{\tilde{B}}(y)=\bigvee \limits_{0\leqslant\alpha\leqslant1}\alpha C_{B(\alpha)}(y)=\bigvee \limits_{0\leqslant\alpha\leqslant1}\mu_{[B_L(\alpha)_{\alpha},B_U(\alpha)_{\alpha}]}(y)
\end{align}
where, $\alpha B(\alpha)$ is a fuzzy set whose membership function is given by
\begin{center}$\mu_{\alpha B(\alpha)}(y) =\begin {cases} \alpha & \mbox {if~~} y\in B(\alpha) \\ 0 & \mbox {if~~}$ otherwise $\end{cases}$
\end{center}
and $C_{B(\alpha)}(y)$ is called characteristic function of $B(\alpha)$ whose membership function is defined as
\begin{center}$C_{B(\alpha)}(y)=\begin {cases} 1 & \mbox {if~~} y\in B(\alpha) \\ 0 & \mbox {if~~} y\notin B(\alpha) \end{cases}$
\end{center}
\textbf{Property 1.} For any a, b, c, d, k $\in \mathbb{R}$ and $a<b, ~c<d$, the interval operations are defined as follows:
\begin{enumerate}[label=(\roman*)]
\item $[a,b]+[c,d]=[a+c,b+d]$
\item $[a,b]-[c,d]=[a-d,b-c]$
\item $k(.)[a,b]=\begin {cases} [ka,kb] & \mbox {if~~} k>0 \\ [kb,ka] & \mbox {if~~} k<0 \end{cases}$\\

Furthermore, for $a>0$ and $c>0$,
\item $[a,b](.)[c,d]=[ac,bd]$
\item $[a,b](\div)[c,d]=[\frac{a}{d}.\frac{b}{c}]$
\end{enumerate}
\subsection{Signed distance (as in Yao and Wu\cite{yaowu})}

For any $a\in\mathbb{R}$, the signed distance from a to 0 is defined as $d_0(a,0)=a$. If $a>0$, then the distance from a to 0 is $a=d_0(a,0)$. If $a<0$, then the distance from a to 0 is $-a=-d_0(a,0)$. Hence, $d_0(a,0)=a$ is known as the signed distance from a to 0.\\

Let $\Omega$ denote the family of all fuzzy sets defined on $\mathbb{R}$. For $\tilde{B}\in \Omega$ with the $\alpha$-cut $B(\alpha)=[B_L(\alpha),B_U(\alpha)],~\alpha \in[0,1]$ where both $B_L(\alpha)$ and $B_U(\alpha)$ are continuous functions on $\alpha\in[0,1]$, the signed distance of the two end-points $B_L(\alpha)$ and $B_U(\alpha)$ of this $\alpha$-cut to the origin 0 is $d_0(B_L(\alpha),0)=B_L(\alpha)$ and $d_0(B_U(\alpha),0)=B_U(\alpha)$ respectively (by the definition stated above). Their average $\frac{B_L(\alpha)+B_U(\alpha)}{2}$ is taken as the signed distance of the $\alpha$-cut $[B_L(\alpha)+B_U(\alpha)]$ to 0, i.e., the signed distance of the interval $[B_L(\alpha),B_U(\alpha)]$ to 0 is defined as
 \begin{center}
 $d_0([B_L(\alpha),B_U(\alpha)],0)=\frac{d_0(B_L(\alpha),0)+d_0(B_U(\alpha),0)}{2}=\frac{B_L(\alpha)+B_U(\alpha)}{2}$
\end{center} 
Since the crisp interval $[B_L(\alpha),B_U(\alpha)]$ has a one-to-one correspondence with the $\alpha$-level fuzzy interval $[B_L(\alpha)_{\alpha},B_U(\alpha)_{\alpha}]$, for every $\alpha \in [0,1]$, we have
\begin{align}
 [B_L(\alpha),B_U(\alpha)]\leftrightarrow  [B_L(\alpha)_{\alpha},B_U(\alpha)_{\alpha}]
 \end{align}
 Also, the real number 0 maps to fuzzy point $\tilde{0}$. Hence, the signed distance of $[B_L(\alpha)_{\alpha},B_U(\alpha)_{\alpha}]$ to $\tilde{0}$ is defined as
\begin{align}
d([B_L(\alpha)_{\alpha},B_U(\alpha)_{\alpha}],\tilde{0})=d([B_L(\alpha)_{\alpha},B_U(\alpha)_{\alpha}],0)=\frac{B_L(\alpha)+B_U(\alpha)}{2}
\end{align}
Using integration, the mean value of the signed distance is obtained as
\begin{align}
\int_{0}^{1}d([B_L(\alpha)_{\alpha},B_U(\alpha)_{\alpha}],\tilde{0})d\alpha=\frac{1}{2}\int_{0}^{1}\{B_L(\alpha)+B_U(\alpha)\}d\alpha
\end{align}
From eqns. (6) and (2), we get
\begin{align}
d(\tilde{B},\tilde{0})=\int_{0}^{1}d([B_L(\alpha)_{\alpha},B_U(\alpha)_{\alpha}],\tilde{0})d\alpha=\frac{1}{2}\int_{0}^{1}\{B_L(\alpha)+B_U(\alpha)\}d\alpha
\end{align}

and accordingly, the following results can also be obtained:\\

\textbf{Result 1.} If $\tilde{B}\in\Omega$ is a triangular fuzzy number such that 
$\tilde{B}=(\beta_1,\beta_2,\beta_3)$, then the $\alpha$-cut of $\tilde{B}$ is $B(\alpha)=[B_L(\alpha)+B_U(\alpha)], ~ \alpha \in [0,1]$, where $B_L(\alpha) = \beta_1+(\beta_2-\beta_1)\alpha$, $~B_R(\alpha)= \beta_3-(\beta_3-\beta_2)\alpha]$, then the signed distance from $\tilde{B}$ to $\tilde{0}$ is $d(\tilde{B},\tilde{0})=\frac{1}{4}(\beta_1+2\beta_2+\beta_3)$.\\

\textbf{Result 2.} For two fuzzy sets$\tilde{B},\tilde{D}\in\Omega$, where $\tilde{B}=\bigcup\limits_{0\leqslant\alpha\leqslant1}[B_L(\alpha)_{\alpha},B_U(\alpha)_{\alpha}]$ and $\tilde{D}=\bigcup\limits_{0\leqslant\alpha\leqslant1}[D_L(\alpha)_{\alpha},D_U(\alpha)_{\alpha}]$, from Property 1 and eqn. (5),we have
\begin{enumerate}
\item $\tilde{B}(+)\tilde{D}=\bigcup\limits_{0\leqslant\alpha\leqslant1}[(B_L(\alpha)+D_L(\alpha))_{\alpha},(B_U(\alpha)+D_U(\alpha))_{\alpha}$
\item $\tilde{B}(-)\tilde{D}=\bigcup\limits_{0\leqslant\alpha\leqslant1}[(B_L(\alpha)-D_L(\alpha))_{\alpha},(B_U(\alpha)-D_U(\alpha))_{\alpha} ]$\\
\item $k(.)\tilde{B}=\begin {cases} \bigcup\limits_{0\leqslant\alpha\leqslant1}[(kB_L(\alpha))_{\alpha},(kB_U(\alpha))_{\alpha}] & \mbox {if~~} k>0\\ \bigcup\limits_{0\leqslant\alpha\leqslant1}[(kB_U(\alpha))_{\alpha},(kB_L(\alpha))_{\alpha}]  & \mbox {if~~} k<0 \\
~~~~~~~~~~~~~~~~~\tilde{0} & \mbox {if~~} k=0\end{cases}~~~~ ,k\in \mathbb{R}$\\
\end{enumerate}

\textbf{Result 3.} For two fuzzy sets $\tilde{B}, \tilde{D}\in\Omega$, and $k\in\mathbb{R}$, from Result 2 and the definition of signed distance, we have
\begin{enumerate}
\item $d(\tilde{B}(+)\tilde{D},\tilde{0})=d(\tilde{B},\tilde{0})+d(\tilde{D},\tilde{0})$
\item $d(\tilde{B}(-)\tilde{D},\tilde{0})=d(\tilde{B},\tilde{0})-d(\tilde{D},\tilde{0})$
\item $d(\tilde{k}(.)\tilde{D},\tilde{0})=kd(\tilde{D},\tilde{0})$
\end{enumerate}

\section{Fuzzy leakage inventory model}
In this fuzzy model, we consider the demand rate $\phi$ and leakage rate $\psi$ of Tomba and Geeta's model to be imprecise in nature. So, we fuzzify them to triangular fuzzy numbers $\tilde{\phi}$ and $\tilde{\psi}$ respectively, defined as $\tilde{\phi} = (\phi-\Delta_{1},\phi,\phi+\Delta_2)$ and $\tilde{\psi} = (\psi-\Delta_3,\psi,\psi+\Delta_4)$, where $0<\Delta_1,\Delta_2<\phi$ and $0<\Delta_3,\Delta_4<\psi$. $\Delta_1,\Delta_2,\Delta_3,\Delta_4$ are the variables to be determined by the decision makers.\\
Then, incorporating the fuzzy demand rate and fuzzy leakage rate in eqn (1),the fuzzy total cost is given by\\
\begin{align}
\tilde{TC} = \frac{qh}{2}+\frac{s\tilde{\phi}}{q}+\frac{qh\tilde{\psi}}{2(\tilde{\phi}+\tilde{\psi})}
\end{align}
Using Result 2, the fuzzy total cost is defuzzified as \\
\begin{align}
Z(q) = d(\tilde{TC},0) = \frac{qh}{2}+\frac{s}{q}d(\tilde{\phi},0)+\frac{qh}{2}d\bigg(\frac{\tilde{\psi}}{\tilde{\phi}+\tilde{\psi}},0\bigg)
\end{align}
Then, \begin{align}
d(\tilde{\phi},\tilde{0}) = \frac{1}{4}[(\phi-\Delta_1) + 2\phi + (\phi+\Delta_2)] = \phi + \frac{(\Delta_2-\Delta_1)}{4}
\end{align}
Let $\tilde{\eta}=\tilde{\phi}+\tilde{\psi}\\
~~~~~~~~= (\phi-\Delta_{1},\phi,\phi+\Delta_2) +(\psi-\Delta_3,\psi,\psi+\Delta_4)\\
~~~~~~~~= ((\phi+\psi)-(\Delta_1+\Delta_3),\phi+\psi,(\phi+\psi)+(\Delta_2+\Delta_4))\\
~~~~~~~~= (\eta-\Delta_5,\eta,\eta_+\Delta_6)$\\
where
\begin{align}
\Delta_5 &=\Delta_1+\Delta_3\\
\Delta_6 &=\Delta_2+\Delta_4
\end{align}
For $(0\leqslant \alpha \leqslant 1)$,\\
Left end-point of $\alpha$-cut of $\tilde{\psi}$ is $\psi_L(\alpha)=(\psi-\Delta_3)+\{\psi-(\psi-\Delta_3)\}\alpha = (\psi-\Delta_3)+\Delta_3\alpha$\\
Right end-point of $\alpha$-cut of $\tilde{\psi}$ is $\psi_U(\alpha)=(\psi+\Delta_4)-\{(\psi+\Delta_4)-\psi\}\alpha = (\psi+\Delta_4)-\Delta_4\alpha$\\
Left end-point of $\alpha$-cut of $\tilde{\eta}$ is $\eta_L(\alpha) = (\eta-\Delta_5)+\{\eta-(\eta-\Delta_5)\}\alpha = (\eta-\Delta_5)+\Delta_5\alpha$\\
Right end-point of $\alpha$-cut of $\tilde{\eta}$ is $\eta_U(\alpha) = (\eta+\Delta_6)-\{(\eta+\Delta_6)-\eta\}\alpha = (\eta+\Delta_6)-\Delta_6\alpha$\\

Since $0< \eta_L(\alpha)<\eta_U(\alpha)$, from Property 1, the left and right end-points of $\alpha$-cut of $\frac{\tilde{\psi}}{\tilde{\phi}+\tilde{\psi}}$ are \\
\begin{align}
\bigg(\frac{\tilde{\psi}}{\tilde{\phi}+\tilde{\psi}}\bigg)_L (\alpha) = \bigg(\frac{\tilde{\psi}}{\tilde{\eta}}\bigg)_L(\alpha)=\frac{\psi_L(\alpha)}{\eta_U(\alpha)}=\frac{(\psi-\Delta_3)+\Delta_3\alpha}{(\eta+\Delta_6)-\Delta_6\alpha} 
\end{align}
\begin{align}
\bigg(\frac{\tilde{\psi}}{\tilde{\phi}+\tilde{\psi}}\bigg)_U(\alpha) = \bigg(\frac{\tilde{\psi}}{\tilde{\eta}}\bigg)_U(\alpha)=\frac{\psi_U(\alpha)}{\eta_L(\alpha)}=\frac{(\psi+\Delta_4)-\Delta_4\alpha}{(\eta-\Delta_5)+\Delta_5\alpha} 
\end{align}
respectively. Then, the signed distance of $\frac{\tilde{\psi}}{\tilde{\alpha}+\tilde{\psi}}$ to ${\tilde 0}$ is \\
\begin{align}
d\bigg(\frac{\tilde{\psi}}{\tilde{\alpha}+\tilde{\psi}},\tilde{0}\bigg) 
&=d\bigg(\frac{\tilde{\psi}}{\tilde{\eta}},\tilde{0}\bigg) \notag\\
&=\frac{1}{2}\int_{0}^{1}\bigg[\bigg(\frac{\psi}{\eta}\bigg)_L(\alpha)+\bigg(\frac{\psi}{\eta}\bigg)_U(\alpha)\bigg]d\alpha\notag \\
&=\frac{1}{2}\int_{0}^{1}\bigg[\frac{(\psi-\Delta_3)+\Delta_3\alpha}{(\eta+\Delta_6)-\Delta_6\alpha} + \frac{(\psi+\Delta_4)-\Delta_4\alpha}{(\eta-\Delta_5)+\Delta_5\alpha}\bigg]d\alpha \notag \\
&= \frac{1}{2}\bigg\{\bigg(\frac{\psi\Delta_6+\eta\Delta_3}{\Delta_6^2}\bigg)\log\bigg(\frac{\eta+\Delta_6}{\eta}\bigg)+\bigg(\frac{\psi\Delta_5+\eta\Delta_4}{\Delta_5^2}\bigg)\log\bigg(\frac{\eta}{\eta-\Delta_5}\bigg)-\frac{\Delta_3}{\Delta_6}-\frac{\Delta_4}{\Delta_5}\bigg\}
\end{align}
which is positive since $\bigg(\frac{\tilde{\psi}}{\tilde{\eta}}\bigg)_L(\alpha)>0$ and $\bigg(\frac{\tilde{\psi}}{\tilde{\eta}}\bigg)_U(\alpha)>0$ are continuous functions on $0\leqslant \alpha \leqslant 1$ and hence the above definite integral is positive.\\
Using the values of $d(\tilde{\phi},\tilde{0})$ and  $d\bigg(\frac{\tilde{\psi}}{\tilde{\alpha}+\tilde{\psi}},\tilde{0}\bigg)$ from eqns. (10) and (15) in eqn. (9), we get\\
\begin{align}
Z(q) = d(\tilde{TC},\tilde{0})=\frac{qh}{2}+\frac{s}{q}\delta+\frac{qh}{2}\zeta
\end{align}
where \begin{align} 
\delta &= d(\tilde{\phi},\tilde{0})=\phi + \frac{(\Delta_2-\Delta_1)}{4}\\
\zeta &= d\bigg(\frac{\tilde{\psi}}{\tilde{\phi}+\tilde{\psi}},\tilde{0}\bigg)= \frac{1}{2}\bigg\{\bigg(\frac{\psi\Delta_6+\eta\Delta_3}{\Delta_6^2}\bigg)\log\bigg(\frac{\eta+\Delta_6}{\eta}\bigg)+\bigg(\frac{\psi\Delta_5+\eta\Delta_4}{\Delta_5^2}\bigg)\log\bigg(\frac{\eta}{\eta-\Delta_5}\bigg)-\frac{\Delta_3}{\Delta_6}-\frac{\Delta_4}{\Delta_5}\bigg\}
\end{align}
Equating the first order partial derivative of Z(q) with (respect to q) to zero, we get 
\begin{center}
$q = \sqrt{\frac{2s\delta}{h(1+\zeta)}}$
\end{center}
Again, taking second order partial derivative of Z(q) with respect to q, we get $\frac{\partial ^2Z}{\partial q^2} = \frac{2s\delta}{q^3}>0$.\\
$\therefore$ q obtained above is the optimal value, i.e., optimal order quantity of the model proposed, i.e., 
\begin{align}
q^* = \sqrt{\frac{2s\delta}{h(1+\zeta)}}    
\end{align}
Also, Z(q) is minimum when $q=q^* = \sqrt{\frac{2s\delta}{h(1+\zeta)}}$. So, from eqn.(16), minimum average total cost
\begin{align}
 Z^*(q)=  d(\tilde{TC},\tilde{0})=\frac{q^*h}{2}+\frac{s}{q^*}\delta+\frac{q^*h}{2}\zeta 
\end{align}
\section{Numerical example and sensitivity analysis}

To illustrate the applicability of the model, the example from Tomba and Geeta\cite{tomba} is considered.\\

The crisp inventory system has the data: 
demand rate $\phi$=600 units per unit per year, holding cost h = \rupee10, set up cost s = \rupee 100 per order, leakage rate $\psi$ = 10 units per unit per year\\
Then, optimal order quantity $q^{*}$= 108 and the minimum cost per year $C_{min} $=\rupee 1104 \\

In our fuzzy model, instead of taking the demand rate and leakage rate as fixed quantity, we assumed them to be imprecise and hence represented by triangular fuzzy numbers $\tilde{\phi} = (600-\Delta_1,600,600+\Delta_2)$ and $\tilde{\psi} = (10-\Delta_3,10,10+\Delta_4)$ respectively in accordance to our fuzzy model. The optimal order lot size $q^*$ and minimum total cost $Z^*$ for various sets of $\Delta_1,\Delta_2,\Delta_3,\Delta_4$ are summarized in Table 1. Relative variation between fuzzy case and crisp case for the optimal lot size $q^*$ and minimum total cost $Z^*$ is also calculated as $Rel_q=\frac{q^*-q^*_c}{q^*_c}\times100\%$ and $Rel_Z=\frac{Z^*-Z^*_c}{Z^*_c}\times100\%$ respectively, where $q^*_c$ and $Z^*_c$ denote the optimal order quantity and minimum total cost in crisp sense respectively.\\
\pagebreak
\setlength{\arrayrulewidth}{0.2mm}
\setlength{\tabcolsep}{10pt}
\renewcommand{\arraystretch}{1.5}
\begin{center}
Table 1. Optimal Solution for the proposed model with fuzzy demand rate and fuzzy leakage rate    
\end{center}
\begin{tabular}{p{0.5cm}p{0.5cm}p{0.3cm} p{0.3cm}p{1cm}p{1.5cm}p{1.5cm}p{1.5cm}p{1.5cm}p{1.5cm}}
\hline
$\Delta_1$ & $\Delta_2$ & $\Delta_3$ & $\Delta_4$ & $d(\tilde{\phi},0)$ & $ d(\frac{\tilde{\psi}}{\tilde{\phi}+\tilde{\psi}},0)$ & $~~~~q^*$ & $~~Z^*$ & $~~{Rel}_q$ & $~~Rel_Z$\\
\hline
100 & 100 & 1 & 1 & 600 & 0.00162 & 109.456 & 1099.76 & 1.3478 & -0.6943\\
~ & ~ & ~ & 2 & 600 & -0.00095 & 109.597 & 1094.92 & 1.4783 & -0.8221\\
~&~& 4 & 3 & 600 & -0.01199 & 110.207 & 1088.86 & 2.0436 & -1.3715\\
~&~&~& 4 & 600 & -0.0144 & 110.342 & 1087.53 & 2.1686 & -1.4921\\
~&~& 7 & 5 & 600 & -0.02498 & 110.939 & 1081.67 & 2.7215 & -2.0223\\
~&~&~& 6& 600 & -0.02725 & 111.068 & 1080.42 & 2.8410 & -2.1362\\
150 & 200 & 1 & 1 & 612.5 &0.00395 & 110.462 & 1108.98 &2.2792 & 0.4514 \\
~&~&~& 2 & 612.5 & 0.0023 & 110.553 & 1108.07 & 2.3636 & 0.3686\\
~&~& 4& 3& 612.5 & -0.00385 & 110.893 & 1104.67 & 2.6790 & 0.0603\\
~&~&~& 4& 612.5 & -0.00544 & 110.982 & 1103.78 & 2.7613 & -0.0199\\
~&~& 7 & 5& 612.5 & -0.01142 & 111.317 & 1100.46 & 3.0717 & -0.3209\\
~&~&~& 6& 612.5 & -0.01296 & 111.404 & 1099.6 & 3.1520 & -0.3986\\
200 & 300 & 1 & 1& 625 & 0.00496 & 111.527 & 1120.8 & 3.2659 & 1.5221\\
~&~&~& 2& 625 & 0.00378 & 111.592 & 1120.15 & 3.3263 & 1.4627\\
~&~& 4&3& 625 & -0.00052 & 111.832 & 1117.75 & 3.5484 & 1.2451\\
~&~&~&  4& 625 & -0.00166 & 111.896 & 1117.11 & 3.6077 & 1.1871\\
~&~&7 & 5 & 625 & -0.00587 & 112.133 & 1114.75 & 3.8268 & 0.9736\\
~&~&~& 6 & 625 & -0.00698 & 112.196 & 1114.12 & 3.8850 & 0.9170\\
\hline
\vspace{0.3cm}
\end{tabular}
From the above Table 1, we observed that 
\begin{enumerate}
\item for fixed $\Delta_1,~\Delta_2$,the demand rate is constant and with the increase in variation of $\Delta_3$ and $\Delta_4$, optimal lot size $q^*$ increases and minimum total cost $Z^*$ decreases. Hence, $Rel_q$ increases and $Rel_Z$ decreases.
\item for fixed $\Delta_3,~\Delta_4$, both the optimal lot size $q^*$ and minimum total cost $Z^*$ increases and hence both $Rel_q$ and $Rel_Z$ increases.
\end{enumerate}

From the discussed example and the sensitivity analysis of the results from Table 1, we cannot exactly ascertain the better results of the model in the two different environments. But still, the one in fuzzy sense will be more reliable as it incorporates real-life situations and conditions with wide range of variability.

\section{Conclusion}

In this paper, we have discussed a leakage inventory model in fuzzy environment and compared it with its traditional model in crisp sense by considering different datas. The main parameters of the model, namely the demand rate and the leakage rate are assumed to be triangular fuzzy numbers. The optimum results are obtained by using signed distance method of defuzzification. It is observed that for different sets of fuzzy demand rate and fuzzy leakage rate, the optimum ordering quantity and minimum total cost are almost more or less with that obtained in crisp environment. Also, uncertainties that are prevalent in real inventory problems are highlighted in this fuzzy model and from the sensitivity analysis by taking different values of $\Delta_1,~\Delta_2,~\Delta_3$ and $\Delta_4$ the variations or the effect of uncertainties the optimum ordering quantity and minimum total cost are analysed.


\begin{thebibliography}{99}
\bibitem{talei}Taleizadeh AA. An EOQ model with partial backordering and advance payments for an evaporating item. International Journal of Production Economics 2014;155(C):185-193
\bibitem{tomba} Tomba I, Geeta O. Some deterministic leakage inventory models. Bulletin of Pure and Applied Sciences 2008;27(2):267-276
\bibitem{yaolee}Yao JS, Lee HM. Fuzzy inventory without backorder for fuzzy order quantity and fuzzy total demand quantity. Computers and Operations Research 2000;27:935-962
\bibitem{lee}Lee HM, Yao JS. Economic production quantity for fuzzy demand quantity and fuzzy production quantity. European Journal of Operational Research 1998;109:203-211
\bibitem{schang}Chang S. Fuzzy production inventory for fuzzy product quantity with triangular fuzzy number. Fuzzy Sets and Systems 1999;107:37-57
\bibitem{yaosu}Yao JS, Su JS. Fuzzy inventory with backorder for fuzzy order quantity. Information Sciences 1996;93:283-319
\bibitem{chang} Chang HC. An application of fuzzy sets theory to the EOQ model with imperfect quality items. Computers \& Operations Research 2004;31:2079-2092
\bibitem{mahata}Mahata GC. Application of fuzzy sets theory in an EOQ model for items with imperfect quality and shortage backordering. International Journal of Services and Operations Management 2013;14(4):466-490
\bibitem{jaggi}Jaggi CK, Anuj S, Mandeep M. A fuzzy inventory model for deteriorating items with initial inspection and allowable shortage under the condition of permissible delay in payment. International Journal of Inventory Control and Management 2012;2(2):167-200
\bibitem{patro}Patro R, Mitali MN, Acharya M. An EOQ model for fuzzy defective rate with allowable proportionate discount. OPSEARCH 2019;56:191-215
\bibitem{chen}Chen LH, Ouyang LY. Fuzzy inventory model for deteriorating items with permissible delay in payment. Applied Mathematics and Computation 2006;182:711-726 
\bibitem{changetal}Chang HC, Yao JS, Ouyang LY. Fuzzy misture inventory model with variable lead-time based on probabilistic fuzzy set and triangular fuzzy number. Mathematical and Computer Modelling 2004;39:287-304
\bibitem{yaowu} Yao JS, Wu K. Ranking fuzzy numbers based on decomposition principle and signed distance. Fuzzy Sets and Systems 2000;116:275-288

\end{thebibliography}
\end{document}